\documentclass[10pt, reqno]{amsart}
\usepackage{amssymb,latexsym}
\usepackage{eucal}

\newtheorem{lemma}{Lemma}

\newtheorem*{Th}{Theorem}

\newcommand{\la}{{\langle}}
\newcommand{\ra}{{\rangle}}
\newcommand{\e}{{\varepsilon}}

\newcommand{\R}{\mathcal{R}}  \newcommand{\Ss}{\mathcal{S}}

\newcommand{\D}{\Delta}
\newcommand{\A}{\mathcal{A}}

\newcommand{\X}{\mathcal{X}}

\newcommand{\p}{\partial}

\newcommand{\C}{\mathcal{C}}

\newcommand{\ph}{\varphi}
\newcommand{\MM}{\mathfrak{M}}
\newcommand{\wtl}{\widetilde}

\begin{document}

\title[Non-hopfian relatively free groups]
{Non-hopfian relatively free groups}
\author{S.V. Ivanov}
\address{Department of Mathematics\\
University of Illinois \\
Urbana,  IL 61801, U.S.A.}
\email{ivanov@math.uiuc.edu}
\thanks{The first author is supported in part by NSF grant
DMS 00-99612}

\author{A.M. Storozhev}
\address{Australian Mathematics Trust\\
University of Canberra \\
Belconnen, ACT 2601, Australia}
\email{andreis@amt.canberra.edu.au }


\begin{abstract}
To solve  problems of  Gilbert Baumslag and Hanna Neumann, posed in the
1960's, we construct a nontrivial variety of groups all of whose noncyclic
free groups are non-hopfian.
\end{abstract}
\maketitle

\section{Introduction}

Recall that a group $G$ is called {\em hopfian} if every
epimorphism $G \to G$ is an automorphism. Let $F_m$ be a free
group of finite rank $m >1$, $N$ be a normal subgroup of $F_m$ and
$V(N)$ a verbal subgroup of $N$ defined by a set of words $V$. In
\cite{B63}, Baumslag proved that if  both quotients $F_m / N$, $N
/ V(N)$ are residually finite then the group $F_m / V(N)$ is also
residually finite. In this connection, Baumslag \cite{B63} posed
the following problem. Is $F_m / V(N)$  a hopfian group if the
group $F_m/N$ is hopfian? In particular,  if  $F_m = N$ then
$F_m/N$ is trivial and hopfian and so the Baumslag problem asks
about the hopfian property of relatively free groups $F_m /
V(F_m)$, where $V(F_m)$ is a verbal (or fully invariant) subgroup
of  $F_m $. The problem on the hopfian property of finitely
generated  relatively free groups was independently stated by H.
Neumann \cite[Problem 15]{N67}. Recall that a finitely generated
residually finite group is hopfian (e.g., see Corollary 41.44
\cite{N67}). Hence, if one could show that relatively free groups
are residually finite then H. Neumann's problem would be solved in
the affirmative. However, this is not the case in general and it
follows from Novikov--Adian results \cite{NA68}, \cite{A75} (see
also \cite{O82}, \cite{O89},  \cite{I94}, \cite{L96}, \cite{I98})
on the Burnside problem for odd $n \ge 665$ and
Kostrikin--Zelmanov results \cite{K59}, \cite{K86}, \cite{Z91} on
the restricted Burnside problem for $n = p^k$, where $p$ is prime,
that the free $m$-generator Burnside group $B(m,n)= F_m /F_m^n$ of
exponent $n$ is not residually finite if $m >1$ and $n$ is odd, $n
= p^k >665$ (in fact,  $B(m,n)$ is not residually finite for all
$m >1$ and $n \gg 1$ as follows from results of \cite{K59},
\cite{K86}, \cite{Z91}, \cite{Z92}, \cite{HH56}, \cite{NA68},
\cite{I94} and the classification of finite simple groups).
Whether the group $B(m,n)$ is hopfian (or, more generally, whether
there is a hopfian non-residually finite relatively free group) is
still unknown (the problem on the hopfian property of $B(m,n)$ for
odd $n \gg 1$ is stated in \cite[Problem 11.36(c)]{KN92}).

In this article we construct a variety of groups of exponent 0 all of whose
noncyclic free  groups are non-hopfian providing thereby  negative solutions
to the foregoing problems of Baumslag and H. Neumann. To construct the
identities that define such a  variety of groups, we let $[a,b] =
aba^{-1}b^{-1}$ be the commutator of $a$, $b$ and set
\begin{equation}
\label{v1v2}
v_0(x,y) =x, \quad  v_1(x,y) =[((x^dy^d)^dx^d)^d, x^d]^d y ,
\quad v_2(x,y) =[v_1(x,y)^d, x^d] .
\end{equation}

Now we define the words $w_1(x,y)$, $w_2(x,y)$ by following formulas
\begin{multline}
\label{w1}
w_1(x,y) = x^{\e_1} v_1(x,y)^{n}  x^{\e_2} v_1(x,y)^{n+2}  \dots
 x^{\e_{h/2-2}} v_1(x,y)^{n+h-6} \\
 x^{\e_{h/2-1}} v_1(x,y)^{n+h-4}
 x^{\e_{h/2}} v_1(x,y)^{(n+h-2) + h/2}
x^{\e_1} v_1(x,y)^{-(n+1)} \\ x^{\e_2} v_1(x,y)^{-(n+3)}  \dots
x^{\e_{h/2-1}} v_1(x,y)^{-(n+h-3)}   x^{\e_{h/2}} v_1(x,y)^{-(n+h-1)} ,
\end{multline}
\begin{multline}
\label{w2}
w_2(x,y) = y v_2(x,y)^{n^2+1}  v_1(x,y)^{\e_2} v_2(x,y)^{n^2+2}
v_1(x,y)^{\e_3} v_2(x,y)^{n^2+3}  \dots   \\
\dots v_1(x,y)^{\e_{h-1}} v_2(x,y)^{n^2+h-1}  v_1(x,y)^{\e_{h}}
v_2(x,y)^{n^2+h} ,
\end{multline}
where $h \equiv 0 \pmod{20}$,
\begin{gather*}
\e_{10k+1}=  \e_{10k+2}=  \e_{10k+3}=  \e_{10k+5}=  \e_{10k+6}=  1, \\
\e_{10k+4}=  \e_{10k+7}=  \e_{10k+8}=  \e_{10k+9}=  \e_{10k+10}=  -1
\end{gather*}
for  $k=  0,1,\dots ,h/10-1$, and $h, d, n$ are sufficiently large
positive integers (with $n \gg d \gg h \gg 1$).  Note that the
sums of exponents on $x, y$ in both (\ref{w1}) and (\ref{w2}) are
zeros.

Now we can state our main result (note that this result  was
announced in \cite[Theorem~5]{IO91}; for the sake of simplicity of
proofs we changed the identities $w_1(x,y) \equiv 1$, $w_2(x,y)
\equiv 1$ presented in \cite{IO91}).

\begin{Th}
Let $\MM$ be the variety of groups defined by identities
$w_1(x, y) \equiv 1$, $w_2(x, y) \equiv 1$, where the words
$w_1(x, y)$,  $w_2(x, y)$ are given by formulas $(\ref{w1})$,
$(\ref{w2})$. Then any free group of rank $m > 1$ in the
variety $\MM$ is not hopfian.
\end{Th}

To prove this Theorem, in Sect.~2,  we inductively construct a
group presentation  by means of generators and defining relations.
In Sect.~3, we apply the geometric machinery of graded diagrams,
developed by Ol'shanskii \cite{O85}, \cite{O89}, to study
subpresentations of this presentation  and to prove a number of
technical lemmas. In particular, we  use the notation and
terminology of \cite{O89} and all notions that are not defined in
this paper can be found in \cite{O89}. In Sect.~4, we show that
the group presentation constructed in Sect.~2 defines a free group
of rank $m > 1$ in the variety $\MM$ and this group is
non-hopfian.

\section{Inductive Construction}

As in \cite{O89}, we will use numerical parameters
$$
\alpha \succ \beta  \succ \gamma  \succ \delta  \succ \e \succ \zeta \succ
\eta \succ \iota
$$
and $h = \delta^{-1}$, $d = \eta^{-1}$, $n = \iota^{-1}$ ($h, d, n$ were
already used in (\ref{w1})--(\ref{w2})) and employ the least parameter
principle (LPP) (according to LPP a small positive value for, say,  $\zeta$
is chosen to satisfy all inequalities whose smallest (in terms of the
relation $\succ$) parameter is $\zeta$).

Let $\A = \{ a_1, \dots, a_m \}$  be an alphabet, $m > 1$, and
$F(\A)$ be the free group in $\A$. Elements of $F(\A)$ are
referred to as {\em words in $\A^{\pm 1} = \A \cup \A^{-1}$} or
just words. Denote $G(0) = F(\A)$ and let the set $\R_0$ be
empty. To define the group $G(i)$ by induction on $i \ge 1$,
assume that the group $G(i-1)$ is already constructed by its
presentation
$$
G(i-1) = \la \A \; \| \; R=1, R \in \R_{i-1} \ra .
$$
Let $\X_i$ be a  set of words (in $\A^{\pm 1}$) of length $i$,
called {\em periods of rank $i$}, which is maximal with respect
to the following two properties:
\begin{enumerate}
\item[(A1)] If $A \in \X_i$ then $A$ (that is, the image of $A$
in $G(i-1)$) is not conjugate in $G(i-1)$ to a power of a word
of length $< |A| = i$.

\item[(A2)] If $A$, $B$ are distinct elements of $\X_i$ then
$A$ is not conjugate in $G(i-1)$ to $B$ or $B^{-1}$.
\end{enumerate}

If the images of two words $X$, $Y$ are equal in the group
$G(i-1)$, $i \ge 1$, then we will say that $X$  is {\em equal
in rank} $i-1$ to $Y$ and will write $ X \overset {i-1} = Y$.
Analogously, we will say that two words $X$, $Y$ are {\em
conjugate in rank} $i-1$ if their images are conjugate in the
group $G(i-1)$. As in  \cite{O89}, a word $A$ is called {\em
simple} in rank $i-1$, $i \ge 1$, if $A$ is conjugate in rank
$i-1$ neither to a power $B^k$, where $|B| <|A|$ nor to a power
of period of some rank $\le i-1$. We will also say that two
pairs $(X_1, X_2)$, $(Y_1, Y_2)$ of words are conjugate in rank
$i-1$, $i \ge 1$, if there is a word $W$ such that  $X_1
\overset {i-1} = W Y_1 W^{-1}$ and $X_2 \overset {i-1} = W Y_2
W^{-1}$.

Consider the set of all possible pairs $(X, Y)$ of words in
$\A^{\pm 1}$ and pick $z^* \in \{1, 2\}$. This set is
partitioned by equivalence $z^*$-classes $\C_\ell(z^*)$, $\ell
=1,2, \dots$, of the equivalence relation $\overset { z^*}
\sim$ defined as follows: $(X_1, Y_1) \overset { z^*} \sim
(X_2, Y_2)$ if and only if the pairs $( v_{z^*}(X_1, Y_1),
w_{z^*}(X_1, Y_1) )$ and $( v_{z^*}(X_2, Y_2), w_{z^*}(X_2,
Y_2))$  are conjugate in rank $i-1$. It is convenient to
enumerate (in some way)
$$
\C_{A^f, 1}(z^*) , \C_{A^f, 2}(z^*), \dots
$$
all $z^*$-classes of pairs $(X, Y)$ such that $w_{z^*}(X, Y)
\overset {i-1} \neq 1$ and  $v_{z^*}(X, Y)$ is conjugate in
rank $i-1$ to some power $A^f$, where $A \in \X_i$ and $f$  are
fixed.

It follows from definitions that every class $\C_{A^f, j}(z^*)$
contains a pair
$$
(X_{A^f, j, z^*}, \bar Y_{A^f, j, z^*} )
$$
with the following properties. The word $X_{A^f, j, z^*}$ is
graphically (that is, letter-by-letter) equal to a power of $
B_{A^f, j, z^*}$, where $ B_{A^f, j, z^*}$ is simple in rank
$i-1$ or a period of rank $\le i-1$; $\bar Y_{A^f, j, z^*}
\equiv Z_{A^f, j, z^*} Y_{A^f, j, z^*} Z_{A^f, j, z^*}^{-1}$,
where the symbol \ '$\equiv$' \  means the graphical equality,
$Y_{A^f, j, z^*}$ is graphically equal to a power of $C_{A^f,
j, z^*}$, where $C_{A^f, j, z^*}$ is simple in rank $i-1$ or a
period of rank $\le i-1$. We can also assume that  if $D_1 \in
\{ A, B_{A^f, j, z^*}, C_{A^f, j, z^*} \}$  is conjugate in
rank $i-1$ to $D_2^{\pm 1}$, where $D_2 \in \{ A, B_{A^f, j,
z^*}, C_{A^f, j, z^*} \}$, then $D_1 \equiv D_2$. Finally, the
word $Z_{A^f, j, z^*}$ is picked for fixed $X_{A^f, j, z^*}$,
$Y_{A^f, j, z^*}$ so that the length $|Z_{A^f, j, z^*}|$ is
minimal (and the pair $(X_{A^f, j, z^*}, Z_{A^f, j, z^*}
Y_{A^f, j, z^*}  Z_{A^f, j, z^*}^{-1})$ belongs to $\C_{A^f,
j}(z^*)$). Similar to \cite{O85}, \cite{O89}, \cite{S94}, the
triple $(X_{A^f, j, z^*}, Y_{A^f, j, z^*}, Z_{A^f, j, z^*} )$
is called an $(A^f, j, z^*)$-{\em triple} corresponding to the
class $\C_{A^f, j}(z^*)$ (in rank $i-1$).

Now for every class $\C_{A^f, j}(z^*)$  we pick a corresponding
$(A^f, j, z^*)$-{triple}
$$
(X_{A^f, j, z^*},  Y_{A^f, j, z^*}, Z_{A^f, j, z^*} )
$$
in rank $i-1$ and construct  a defining word $R_{A^f, j, z^*}$
of rank $i$ as follows. Pick a word $W_{A^f, j, z^*}$ of
minimal length so that
$$
v_{z^*}( X_{A^f, j, z^*},  \bar Y_{A^f, j, z^*} ) \overset {i-1} =
W_{A^f, j, z^*} A^{f} W_{A^f, j, z^*}^{-1} .
$$
Let   $T_{A^f, j, z^*}$,   $U_{A^f, j, z^*}$ be words of
minimal length such that
\begin{gather*}
T_{A^f, j, z^*}  \overset {i-1} =
W_{A^f, j, z^*}^{-1}  v_{z^*-1}( X_{A^f, j, z^*},  \bar Y_{A^f, j, z^*} )
W_{A^f, j, z^*} \\
U_{A^f, j, z^*}  \overset {i-1} = W_{A^f, j, z^*}^{-1}  \bar Y_{A^f, j, z^*}
W_{A^f, j, z^*}
\end{gather*}
(recall that $v_0(x,y) = x$, see (\ref{v1v2})).

 If $z^* =1$ then, in accordance with (\ref{w1}), we set
\begin{multline}
\label{R1}
R_{A^f, j, 1}  =   T_{A^f, j, 1}^{\e_1}  A^{ n f } T_{A^f, j,
1}^{\e_2} A^{ (n+2) f }  \dots T_{A^f, j, 1}^{\e_{h/2-2}}  A^{ (n+h-6) f }
T_{A^f, j, 1}^{\e_{h/2-1}}  A^{ (n+h-4) f }   \\
T_{A^f, j, 1}^{\e_{h/2}} A^{ (n+h-2 + h/2) f }
T_{A^f, j, 1}^{\e_1}  A^{ -(n+1) f }
T_{A^f, j, 1}^{\e_2}  A^{ -(n+3) f } \dots      \\
\dots  T_{A^f, j, 1}^{\e_{h/2-1}} A^{ -(n+h-3) f }
T_{A^f, j, 1}^{\e_{h/2}} A^{ -(n+h-1)
f } ,
\end{multline}
where $\e_1, \dots, \e_{h/2}$, $h, n$ are defined as in (\ref{w1}).

If $z^* =2$ then, in accordance with (\ref{w2}), we put
\begin{multline}
\label{R2}
R_{A^f, j, 2}  =   U_{A^f, j, 2}  A^{ (n^2+1) f } T_{A^f, j,
2}^{\e_2}  A^{ (n^2+2) f } T_{A^f, j, 2}^{\e_3}  A^{ (n^2+3) f }  \dots \\
\dots T_{A^f, j, 2}^{\e_{h-1}}  A^{ (n^2+h-1) f } T_{A^f, j, 2}^{\e_{h}}  A^{
(n^2+h) f }   ,
\end{multline}
where $\e_2, \dots, \e_h$, $h, n$ are defined as in (\ref{w2}).

It follows from definitions that the word $R_{A^f, j, z^*}$ is conjugate in
rank $i-1$ (by the word $W_{A^f, j, z^*}^{-1}$) to $w_{z^*}( X_{A^f, j, z^*},
\bar Y_{A^f, j, z^*}) \overset {i-1} \neq 1$.

The set $\Ss_i$ of defining words of rank $i$ consists of all possible words
$R_{A^f, j, z^*}$  given by (\ref{R1})--(\ref{R2}) (over all equivalence
classes $\C_{A^f, j}(z^*)$, $A \in \X_i$). Finally, we put $\R_i = \R_{i-1}
\cup \Ss_i$ and set
\begin{gather}
\label{Gi} G(i) = \la \A \; \| \; R=1, R \in  \R_i \ra .
\end{gather}
 The inductive definition of  groups $G(i)$, $i \ge 0$, is
 now complete and we can consider the limit group $G(\infty)$ given by
 defining words of all ranks $j =1,2, \dots$
\begin{gather}
\label{G8} G(\infty) = \la \A \; \| \; R=1, R \in
\cup_{j=0}^\infty \R_j \ra .
\end{gather}

\section{Several Lemmas}

We will prove (Lemma \ref{L11}) that the group $G(\infty)$,
defined by \eqref{G8}, is the free group of the variety $\MM$ in
the alphabet $\A$, that is, $G(\infty)$ is naturally isomorphic to
the quotient $F(\A) / W_{1,2}(F(\A))$, where $W_{1,2}(F(\A))$ is
the verbal subgroup of $F(\A)$, defined by the set $W_{1,2} = \{
w_1(x,y), w_2(x,y) \}$. We will also show that $G(\infty)$ is not
hopfian (Lemma \ref{L12}). But first we need to study the
presentation (\ref{Gi}) of $G(i)$ and establish a number of
technical lemmas.

As in Sects. 29--30 \cite{O89}, the following Lemmas
\ref{L1}--\ref{L10} are proved by induction on $i \ge 0$ (whose
base for $i =0$ is trivial).
\begin{lemma} 
\label{L1} The presentation $(\ref{Gi})$ of $G(i)$ satisfies the
condition $R$ of \cite[Sect. 25]{O89}.
\end{lemma}

\begin{proof} This proof is quite similar to the proof of Lemma 29.4
\cite{O89}. Inductive references to Lemmas 30.3, 30.4, 30.5 \cite{O89} (in
rank $i-1$) are replaced by references to Lemma \ref{L6}. Note that, by Lemma
\ref{L6}, we have that
$$
n^2 >  100\zeta^{-1}(n+h) \ge  |f|(n +h)
$$
(LPP: $\delta = h^{-1} \succ \zeta \succ \iota = n^{-1}$) which implies that,
repeating the arguments of Lemma 29.3  \cite{O89}, we can conclude that  the
defining relators $R$, $R'$ correspond to the same value of $z^*$. Therefore,
Lemma \ref{L10} enables us to finish the proof of the analog of Lemma 29.3 as
in \cite{O89} (note that we need Lemma \ref{L10} only when $z^* =2$).
\end{proof}

\begin{lemma}
\label{L2A} Suppose that $X_1$, $X_2$ are some words and $k \ne 0$ is an
integer. Then

$(a)$ If $[ X_1^k, X_2^k ] \overset i = 1$, then $[ X_1, X_2 ] \overset i =
1$;

$(b)$ If  $X_1^{k} \overset i = X_2^k$, then $X_1 \overset i = X_2$. In
particular, the group $G(i)$, defined by $(\ref{Gi})$, is torsion-free.
\end{lemma}

\begin{proof}  (a)  It follows from definitions and Lemma 25.2 \cite{O89} that
the group $G(i)$, defined by (\ref{Gi}),  is torsion-free. Hence, we can
apply Lemma 25.12 \cite{O89} to equality  $[ X_1^k, X_2^k ] \overset i = 1$
and obtain that $[ X_1, X_2 ]\overset i = 1$, as required.

(b) The equality $X_1^{k} \overset i = X_2^k$ implies that $[ X_1^k, X_2^k ]
\overset i = 1$ and, by part (a), we have $[ X_1, X_2 ]\overset i = 1$.
Hence, $(X_1 X_2^{-1})^{k} \overset i = 1$ and, as above, it follows from
definitions and Lemma 25.2 \cite{O89} that $X_1 X_2^{-1}  \overset i = 1$.
Lemma \ref{L2A} is proved.
\end{proof}

Recall that a subgroup $H$ of a group $G$ is called {\em antinormal} if for
every $g \in G$ the inequality $g H g ^{-1} \cap H \neq \{ 1\}$ implies that
$g \in H$.

\begin{lemma}
\label{L2B} $(a)$  Every word is conjugate in rank $i$ to a power of either a
simple in rank $i$ word or a period of rank $\le i$.

 $(b)$  Suppose that each of $A$, $B$ is either a simple
 in rank $i$ word or a period of rank $\le i$ and $A^k$ is conjugate
 in rank $i$ to $B^\ell$, $\ell \neq 0$. Then $A$ is conjugate in rank $i$
 to $B$ or to $B^{-1}$.

 $(c)$ Let $A$ be a simple in rank $i$ word or a period of rank $\le i$.
 Then  the cyclic  subgroup $\la A\ra$, generated by $A$, of the group
 $G(i)$, defined by $(\ref{Gi})$, is antinormal.
\end{lemma}

\begin{proof}
Part (a) follows from definitions (see also Lemma 18.1 \cite{O89}). Since the
group $G(i)$ is torsion-free by Lemma \ref{L2A},  we can argue as in the
proof of Lemma 25.17 \cite{O89} (see also the proof of Theorem 19.4
\cite{O89}) to prove part (b) and obtain that an equality of the form $ZA^k
Z^{-1} \overset i = A^\ell$, $\ell \neq 0$, implies that $Z \in \la A \ra$,
as required in part (c).  Lemma \ref{L2B} is proved.
\end{proof}

In addition to $A$-, $B$-, $\dots$, $H$-maps which are introduced and
investigated in \cite{O89}, we will need $I$-maps (cf. \cite{S94}) defined as
follows. A $B$-map $\D$ is called an {\em $I$-map} if following properties
(I1)--(I5) hold.

\begin{enumerate}
\item[(I1)] $\D$ is a map on a sphere punctured at least once and at most
thrice.

\item[(I2)] The cyclic sections of the boundary $\p \D$ of $\D$ are products
of sections (some of which or all can be cyclic) of two types: long sections
and short sections.

\item[(I3)] If $s$ is a long section of $\p \D$ then $s$ is smooth of rank
$r(s)$ and $|s| > 10 \zeta^{-1} r(s)$.

\item[(I4)] If $t$ is a short section of $\p \D$ then $|t| < \zeta |s_0|$,
where $s_0$ is a long section of minimal length.

\item[(I5)] If $L(\p\D)$ and $S(\p\D)$ are the numbers of long and short
sections of $\p \D$, respectively, then  $S(\p\D) \le 2 L(\p\D)$ and $1 \le
L(\p\D) \le 10$.
\end{enumerate}

\begin{lemma} \label{L2I}
Suppose that $\D$ is an $I$-map. Then there is a system of  pairwise disjoint
regular contiguity submaps of long sections to long sections in $\D$ such
that no two distinct contiguity submaps of a long section $s_1$ to a long
section $s_2$ are contained in any larger contiguity submap of $s_1$ to $s_2$
and the sum of contiguity arcs of contiguity submaps of the system is greater
than $(1-\alpha^{1/4})L_s$, where $L_s$ is the sum of lengths of all long
sections of $\p \D$.
\end{lemma}
\begin{proof}
Without loss of generality, we can assume that every short section $t$ of $\p
\D$ is geodesic in $\D$ (and that if $t$ is cyclic then $t$ is cyclically
geodesic; note that we can always replace $t$ by a homotopic to $t$ in $\D$
geodesic path).

This proof is analogous to the proof of Lemma 24.6 \cite{O89} (see also
Lemmas 23.15, 24.2 \cite{O89} on $C$- and $D$-maps). Repeating arguments of
the proof of Lemma 24.6 \cite{O89}, we can establish similar estimates for an
$I$-map $\D$. Note that we need to make straightforward corrections of the
number of distinguished contiguity submaps between $p$ and $q$, where $p$,
$q$ are sections of $\p \D$ or $\p \Pi$ and $\Pi$ is a 2-cell of $\D$, and of
the number of distinguished contiguity submaps between sections of $\p \D$.
As in the proof of Lemma 24.6 \cite{O89}, we obtain the estimate $M < \alpha
\nu(\D)$, where $M$ is the sum of weights of all inner edges of $\D$ and
$\nu(\D)$ is the total weight of $\D$.

Now we can argue as in the proof of Lemma 23.15 \cite{O89} to
derive that the sum of lengths of outer arcs of long sections
of $\p \D$ is greater than $(1-\alpha^{1/4})L_s$.
\end{proof}

\begin{lemma} 
\label{L2C} Suppose that each of $A_1$, $A_2$, $B$ is a simple in rank $i$
word or a period of rank $\le i$, $[ A_1^{k_1} , Z A_2^{k_2} Z^{-1} ]
\overset i \neq 1$ for some word $Z$,  and $B^\ell$ is conjugate in rank $i$
either to $[ A_1^{d k_1} , Z A_2^{d k_2} Z^{-1} ]$ or to $ A_1^{dk_1} Z
A_2^{dk_2} Z^{-1}$. Then $0 < |\ell | \le 100 \zeta^{-1}$ and either $\max (
|A_1^{d k_1}|, |A_2^{d k_2}| ) \le \zeta^{-1} |B^\ell |$ if $B^\ell$ is
conjugate in rank $i$ to $[ A_1^{dk_1} , Z A_2^{dk_2} Z^{-1} ]$ or $\max (
|A_1^{dk_1}|, |A_2^{dk_2}| ) \le \zeta^{-2} |B^\ell |$ if $B^\ell$ is
conjugate in rank $i$ to $ A_1^{dk_1} Z A_2^{dk_2} Z^{-1}$.
\end{lemma}

\begin{proof} Without loss of generality (see also Lemma \ref{L2B}), we can
assume that if $B_1 \in \{ A_1^{\pm 1}, A_2^{\pm 1}, B^{\pm 1} \}$ is
conjugate in rank $i$ to $B_2 \in \{ A_1^{\pm 1}, A_2^{\pm 1}, B^{\pm 1} \}$,
then $B_1 \equiv B_2$.

If $\ell = 0$, that is, either $[ A_1^{dk_1} , Z A_2^{dk_2} Z^{-1} ] \overset
i = 1$ or $ A_1^{dk_1} Z A_2^{dk_2} Z^{-1} \overset i = 1$ then, by Lemma
\ref{L2A}, we have that either $[ A_1^{k_1} , Z A_2^{k_2} Z^{-1} ] \overset i
= 1$ or $ A_1^{k_1} Z A_2^{k_2} Z^{-1} \overset i = 1$, contrary to lemma's
hypothesis   $[ A_1^{k_1} , Z A_2^{k_2} Z^{-1} ] \overset i \neq 1$. Hence,
$\ell \ne 0$.

First assume that
\begin{equation}
\label{L2C:1}
 [ A_1^{dk_1} , Z A_2^{dk_2} Z^{-1} ] \overset i = W_B B^\ell W_B^{-1}
\end{equation}
for some word $W_B$. Then there is a reduced diagram $\D_1$ of rank $i$ on a
thrice punctured sphere the labels of 3 cyclic sections of whose boundary $\p
\D_1$ are $A_1^{d k_1}$, $A_1^{-d k_1}$, $B^\ell$ ($\D_1$ can be constructed
from a simply connected diagram $\D_0$ of rank $i$ for equality (\ref{L2C:1})
by identifying sections of $\p \D_0$ labelled by $Z A_2^{dk_2} Z^{-1}$ and $Z
A_2^{-dk_2} Z^{-1}$, $W_B$ and $W_B^{-1}$). If $|\ell | > 100 \zeta^{-1}$,
then $\D_1$ is a $G$-map (see \cite[Sect. 24.2]{O89}) and, as in the proof of
25.19 \cite{O89}, it follows from Lemma 24.8 \cite{O89} that $B^\ell \overset
i = 1$, contrary to Lemma \ref{L2A} and $\ell \ne 0$. Hence, $|\ell | \le 100
\zeta^{-1}$, as desired.

Suppose  that $|B^\ell | < \zeta | A_1^{d k_1} |$. Then $\D_1$ is an $E$-map
(see \cite[Sect. 24.2]{O89}) and it follows from Lemmas 24.6, 25.10
\cite{O89} that cyclic sections of $\p \D_1$, labelled by $A_1^{d k_1}$,
$A_1^{-d k_1}$, are $A_1$-compatible. Now it is easy to see that $B^\ell
\overset i = 1$, whence $\ell = 0$ by Lemma \ref{L2A}. This contradiction to
$\ell \ne 0$ shows that  $| A_1^{d k_1} | \le \zeta^{-1} |B^\ell|$, as
desired. Analogously, $| A_2^{d k_2} | \le \zeta^{-1} |B^\ell|$ (using
equality (\ref{L2C:1}), we can construct a similar diagram $\D_2$ the labels
of 3 cyclic sections of whose boundary $\p \D_2$ are $A_2^{d k_2}$, $A_2^{-d
k_2}$, $B^\ell$, and then argue as before).

Now assume that
\begin{equation}
\label{L2C:2}
  A_1^{dk_1} Z A_2^{dk_2} Z^{-1} \overset i = W_B B^\ell W_B^{-1}
\end{equation}
for some word $W_B$. Then there is a reduced diagram $\D$ of rank $i$ on a
thrice punctured sphere the labels of 3 cyclic sections of whose boundary $\p
\D$ are $A_1^{d k_1}$, $A_2^{d k_2}$, $B^{-\ell}$ ($\D$ can be constructed
from a simply connected diagram $\D_0$ of rank $i$ for equality (\ref{L2C:2})
by identifying the sections of $\p \D_0$ labelled by $Z$ and  $Z^{-1}$, $W_B$
and $W_B^{-1}$). If $|\ell | > 100 \zeta^{-1}$, then $\D$ is a $G$-map  and,
using Lemmas 24.8, 25.10 \cite{O89}, we can conclude that one of cyclic
sections of $\p \D$ is compatible with another cyclic section of $\p \D$.
Then it follows from Lemmas 3 and 24.9 \cite{O89} that $[ A_1^{dk_1} , Z
A_2^{dk_2} Z^{-1} ] \overset i = 1$ and so, by Lemma \ref{L2A}, $[ A_1^{k_1}
, Z A_2^{k_2} Z^{-1} ] \overset i = 1$, contrary to lemma's hypothesis. This
contradiction shows that $|\ell | \le 100 \zeta^{-1}$.

Suppose that $|B^\ell | < \zeta \min ( |  A_1^{dk_1}|, |A_2^{dk_2}|)$. Then
$\D$ is an $E$-map (see \cite[Sect. 24.2]{O89}) and, it follows from Lemmas
24.6, 25.10 \cite{O89} that  two distinct cyclic sections of $\p \D$ are
compatible. As above, this implies that $[ A_1^{k_1} , Z A_2^{k_2} Z^{-1} ]
\overset i = 1$, contrary to lemma's hypothesis. Hence, it is shown that
$|B^\ell | \ge  \zeta \min ( |  A_1^{dk_1}|, |A_2^{dk_2}|)$. For
definiteness, let $|A_1^{dk_1}| \le |A_2^{dk_2}| $. Then
\begin{equation}
\label{L2C:3} \zeta |A_1^{dk_1}| \le |B^\ell | .
\end{equation}
If $ |A_2^{dk_2}| \le \zeta^{-2}  |B^\ell | $ then our proof is obviously
finished. So we may assume that
\begin{equation}
\label{L2C:4}
 |A_2^{dk_2}| > \zeta^{-2}  |B^\ell |  .
\end{equation}
It follows from inequalities (\ref{L2C:3})--(\ref{L2C:4}) that $|A_1^{dk_1}|,
 |B^\ell | <\zeta |A_2^{dk_2}|$ and we can see that $\D$ is an $I$-map.
 This, however, is impossible by Lemmas \ref{L2I} and 25.10 \cite{O89}.
 This contradiction completes the proof of Lemma  \ref{L2C}.
\end{proof}

Now suppose that $X, Y$ are some words and
\begin{equation}
\label{XY1}
 [X, Y] \overset {i} \neq 1 .
\end{equation}

Conjugating the pair $(X, Y)$ in rank $i$ if necessary, we can assume that $X
\equiv B^{k_B}$,  $Y \equiv Z C^{k_C} Z^{-1}$,  where each of  $B, C$ is
either simple in rank $i$ or a period of some rank $\le i$ and, when
$B^{k_B}$,  $C^{k_C}$ are fixed, the word $Z$ is picked to have minimal
length. Furthermore, consider the following equalities
\begin{gather*}
X^dY^d \overset {i} = W_D D^{k_D} W_D^{-1} ,  \qquad
(X^dY^d )^d X^d \overset {i} = W_E E^{k_E} W_E^{-1} ,  \\
[((X^dY^d )^d X^d)^d,  X^d]  \overset {i} = W_F F^{k_F} W_F^{-1} , \quad
[((X^dY^d )^d X^d)^d,  X^d]^dY  \overset {i} = W_G G^{k_G} W_G^{-1} , \\
[([((X^dY^d )^d X^d)^d,  X^d]^dY)^d,X^d]  \overset {i} = W_H H^{k_H}
W_H^{-1} ,
\end{gather*}
where each of $ D, E, F, G, H$ is either simple in rank $i$ or
a period of some rank $\le i$ and the conjugating words $W_D,
W_E, W_F, W_G, W_H$ are picked (when $D, E, F, G, H$ are fixed)
to have minimal length.

Without loss of generality, we can also assume that if $A_1 \in \{ B, C, D,
E, F, G, H\}$ is conjugate in rank $i$ to $A_2^{\pm 1}$, where $A_2 \in \{ B,
C, D, E, F, G, H\}$, then $A_1 \equiv A_2$.

\begin{lemma} 
\label{L2}

In the foregoing notation, the following estimates hold
\begin{gather}
0 < | k_D | \le 100 \zeta^{-1} ,   \quad  \max(| B^{d k_B} |,  | C^{d k_C} |)
\le
\zeta^{-2} | D^{k_D} | , \label{L2:1} \\
|Z| < 3 \zeta^{-2} | D^{k_D} | , \quad
|W_D| < 5 \zeta^{-2} | D^{k_D} |, \label{L2:2} \\
0 < | k_E | \le 100 \zeta^{-1} , \quad
| D^{d k_D} | \le   \zeta^{-2} | E^{k_E} | , \label{L2:3}   \\
|W_E| < 2\zeta^{-2}| E^{k_E} |, \label{L2:4} \\
0 < | k_F | \le 100 \zeta^{-1} , \quad
| E^{d k_E} | \le   \zeta^{-1} | F^{k_F} | , \label{L2:5}   \\
|W_F| < 2 \zeta^{-1} | F^{k_F} |. \label{L2:6}
\end{gather}
\end{lemma}

\begin{proof} If $k_D = 0$, that is, $X^dY^d  \overset {i} = 1$,  then, by
Lemma \ref{L2A},  $XY \overset {i} = 1$ and, therefore,  $[X, Y ] \overset
{i} = 1$, contrary to inequality (\ref{XY1}). Hence $k_D \neq 0$.

In view of inequality (\ref{XY1}), we can apply Lemma \ref{L2C} to the pair
$X \equiv B^{k_B}$, $Y \equiv Z C^{k_C} Z^{-1}$ which yields that $| k_D |
\le 100 \zeta^{-1}$ and  $| B^{d k_B} |,  | C^{d k_C} |\le \zeta^{-2} |
D^{k_D} |$. Inequalities (\ref{L2:1}) are proved.

In view of equality $B^{d k_B}Z C^{d k_C} Z^{-1} \overset {i} = W_D D^{k_D}
W_D^{-1}$, there is a reduced diagram $\D$ of rank $i$ on a thrice punctured
sphere the labels of three cyclic sections of whose boundary $\p \D$ are
$B^{d k_B}$, $C^{d k_C}$, $D^{-k_D}$. It follows from Lemmas 22.2, 24.9
\cite{O89} that
$$
|Z| < (1+4\gamma) (   | B^{d k_B} |+ | C^{d k_C} | + | D^{k_D} |  ) .
$$
In view of inequalities (\ref{L2:1}), we have
$$
|Z| < (1+4\gamma) ( 2\zeta^{-2}    +1  ) | D^{k_D} | < 3 \zeta^{-2} |
D^{k_D}| ,
$$
as claimed in (\ref{L2:2}).

By estimates (\ref{L2:1}) and already proven inequality $|Z| < 3 \zeta^{-2}
|D^{k_D} |$, we have
$$
|X^dY^d| = | B^{d k_B} |+  | C^{d k_C} | + 2 | Z |  < 8\zeta^{-2} | D^{k_D} |
.
$$
Hence, it follows from Lemmas \ref{L1} and 22.1 \cite{O89} that
$$
|W_D| <  (\gamma + \tfrac 1 2)( | X^dY^d | + | D^{k_D} | ) < 5 \zeta^{-2} |
D^{k_D} |
$$
and inequalities (\ref{L2:2}) are  proved.

Now assume that $k_E=0$. Then $(X^dY^d)^dX^d \overset {i} = 1$ which implies
that
$$
[(X^dY^d)^d,X^d]\overset {i} = 1 .
$$
Hence, by Lemma \ref{L2A}, we have $[X, Y ] \overset {i} = 1$, contrary to
inequality (\ref{XY1}).

Consider the equality
\begin{gather}
\label{WDB}
W_D D^{dk_D} W_D^{-1}B^{d k_B} \overset {i} = W_E E^{k_E}
W_E^{-1} .
\end{gather}

In view of (\ref{L2:1}),
$$
| B|, |C| \le d^{-1} \zeta^{-2}| D^{k_D} | \le d^{-1} \zeta^{-2} \cdot 100
\zeta^{-1}| D |
$$
and so $| B|, |C| < | D |$ (LPP: $\zeta \succ \eta = d^{-1}$). Hence, $ [W_D
D^{dk_D} W_D^{-1}, B^{d k_B} ] \overset {i} \neq 1$ (otherwise, we would have
a contradiction to Lemma \ref{L2B}). This last inequality enables us to apply
Lemma \ref{L2C}  to equality (\ref{WDB}) and conclude that $| k_E | \le 100
\zeta^{-1}$, $| D^{d k_D} | \le \zeta^{-2} | E^{k_E} |$. Inequalities
(\ref{L2:3}) are proved.

By (\ref{L2:1}), (\ref{L2:2}), (\ref{L2:3}), we obtain
$$
|W_D D^{dk_D} W_D^{-1}B^{d k_B}| < (1+11\zeta^{-2}d^{-1})| D^{dk_D} |
<2\zeta^{-2}| E^{k_E} |
$$
(LPP: $\zeta \succ \eta = d^{-1}$). Therefore, it follows from Lemmas
\ref{L1} and 22.1 \cite{O89} that
$$
|W_E| <  (\gamma + \tfrac 1 2)( | W_D D^{dk_D} W_D^{-1}B^{d k_B} | + |
E^{k_E} | )
< 2\zeta^{-2}| E^{k_E} |
$$
as claimed in (\ref{L2:4}).

Next assume that $k_F=0$. Then $[((X^dY^d )^d X^d)^d,  X^d]\overset {i} = 1$.
Hence, by Lemma \ref{L2A}, we obtain $[X, Y ] \overset {i} = 1$, which
contradicts inequality (\ref{XY1}).

Now consider the equality $[ W_E E^{dk_E} W_E^{-1} , B^{d k_B}] \overset {i}
= W_F F^{k_F} W_F^{-1}$. By Lemma \ref{L2C}, $| k_F | \le 100 \zeta^{-1}$, $|
E^{d k_E} | \le \zeta^{-1} | F^{k_F} |$, and estimates (\ref{L2:5}) are
proved.

In view of equality
$[ W_E E^{dk_E} W_E^{-1} , B^{d k_B}] \overset {i} = W_F F^{k_F} W_F^{-1}$,
it follows from Lemma 22.1 \cite{O89} that
$$
| W_F | <   ( \gamma + \tfrac 12 )\cdot 2( 2| W_E | +
 | E^{d k_E} | + | B^{d k_B } |  +  \tfrac 12  |  F^{k_F}  |) .
$$
Hence, by estimates (\ref{L2:1}), (\ref{L2:3}), (\ref{L2:4}), (\ref{L2:5}),
we get
\begin{multline*} | W_F | <   (1 + 2 \gamma ) (  (2 \zeta^{-2} d^{-1} +1
)| E^{d k_E} | + \zeta^{-4} d^{-2}| E^{d k_E} | +   \tfrac 12  |  F^{k_F}  |)
< \\
< (1 + 2 \gamma ) ( ( (2 \zeta^{-2} d^{-1} +1 )  + \zeta^{-4}
d^{-2} ) \zeta^{-1} +  \tfrac 12)  |  F^{k_F}  |) <   2
\zeta^{-1}| F^{k_F} |
\end{multline*} (LPP: $\gamma \succ \zeta \succ
\eta = d^{-1}$), as required.  Lemma \ref{L2} is proved.
\end{proof}

\begin{lemma}
\label{L4}
 In the foregoing notation, the following inequalities hold
\begin{gather}
0 < | k_G | \le 10 \zeta^{-1} , \label{L4:1} \\
| F^{d k_F} |  \le  \zeta^{-1} | G^{k_G} | ,  \label{L4:2}  \\
|W_G| <  \zeta^{-1} | G^{k_G} | . \label{L4:3}
\end{gather}
\end{lemma}

\begin{proof}  If $k_G =0$ then $[((X^dY^d )^d X^d)^d,  X^d]^dY \overset i =
1$ which implies that $F^{dk_F}$ is conjugate in rank $i$ to $C^{-k_C}$. It
follows from definitions and Lemma \ref{L2B} that $F\equiv C$. However, it
follows from Lemma \ref{L2} that
\begin{equation}
|C|\le \zeta^{-5}d^{-3}|F^{k_F}|\le 100\zeta^{-6}d^{-3}|F|<|F|
\label{L4:4.0}
\end{equation}
(LPP: $\zeta \succ \eta = d^{-1}$). Hence $k_G \neq 0$.

By definitions, we have
\begin{equation}
W_F F^{d k_F} W_F^{-1} Z C^{ k_C} Z^{-1} \overset {i} = W_G G^{k_G} W_G^{-1}
. \label{L4:4}
\end{equation}
By Lemma \ref{L2},
\begin{multline}
| W_F^{-1} Z  C^{ k_C} Z^{-1} W_F | <
(4\zeta^{-1} + 6 \zeta^{-5}d^{-2} + \zeta^{-5}d^{-3}) | F^{k_F}  |  < \\
< 5 \zeta^{-1} d^{-1} | F^{d k_F}  | < \zeta  | F^{d k_F}  | \label{L4:5}
\end{multline}
(LPP: $\zeta \succ \eta = d^{-1}$).

If  $| G^{k_G}  | <  \zeta | F^{d k_F}  | $ then a reduced annular diagram of
rank $i$ for conjugacy of words $ F^{d k_F} W_F^{-1} Z C^{ k_C} Z^{-1} W_F$ \
and  \ $G^{k_G}$ (see (\ref{L4:4}))  is an $F$-map (see Sect. 24.2
\cite{O89}) whose existence contradicts Lemma \ref{L1}  and Lemmas 24.7,
25.10 \cite{O89}. Therefore,
$ |  G^{k_G}  | \ge  \zeta | F^{d k_F}  | $
and (\ref{L4:2}) is proved.

In view of equality (\ref{L4:4}), there is a reduced diagram $\D$ of rank $i$
on a thrice punctured sphere the labels of whose cyclic sections  are $F^{d
k_F}$, $C^{k_C}$, $G^{-k_G}$. It follows from Lemma \ref{L2} that
$$
|  C^{k_C} | < \zeta^{-5}d^{-4}|F^{dk_F}| < \zeta^2 | F^{d k_F} |
$$
(LPP: $\zeta \succ \eta = d^{-1}$).  Hence, by (\ref{L4:2}), $| C^{k_C} | <
\zeta \min( | F^{d k_F} | , | G^{k_G}| )$. If $| k_G | >  10 \zeta^{-1}$,
then  $\D$ is an $E$-map (see Sect. 24.2 \cite{O89}) and we can argue as in
the proof of Lemma 25.19 \cite{O89}  to show that $F \equiv G$ and the cyclic
sections of $\D$ labelled by $F^{d k_F}$, $G^{-k_E}$ are $F$-compatible. Then
$C^{k_C}$ is  conjugate in rank $i$ to a power of $F$. This, however, is
impossible by Lemma \ref{L2B} and inequality (\ref{L4:4.0}). Hence, $| k_G |
\le 10 \zeta^{-1}$ and inequalities (\ref{L4:1}) are proved.

By Lemma \ref{L1}, we can apply Lemma  22.1 \cite{O89}  to  a reduced annular
diagram of rank $i$ for conjugacy  of words $W_F F^{d k_F} W_F^{-1} Z C^{
k_C} Z^{-1}$ and $G^{k_G}$ to obtain, using estimates (\ref{L4:2}),
(\ref{L4:5}), that
$$
| W_G | <   ( \gamma + \tfrac 12 )( \zeta+1) (| F^{d k_F} | + | G^{k_G} | )
\le ( \gamma + \tfrac 12 )( \zeta+1) ( \zeta^{-1} +1)| G^{k_G} | < \zeta^{-1}
| G^{k_G} | ,
$$
as claimed in  (\ref{L4:3}). Lemma \ref{L4} is proved.
\end{proof}

\begin{lemma} 
\label{L5}
In the foregoing notation, the following inequalities hold
\begin{gather}
0 < | k_H | \le 100 \zeta^{-1} , \quad
| G^{d k_G} |  \le  \zeta^{-1} | H^{k_H} | , \label{L5:1}  \\
|W_H| <  2   \zeta^{-1}   | H^{k_H} | . \label{L5:2}
\end{gather}
\end{lemma}

\begin{proof}  Assume that  $[([((X^dY^d )^d X^d)^d,  X^d]^dY)^d,X^d]
\overset {i} = 1$. Then, by Lemma \ref{L2A}, we also have $[ [((X^dY^d )^d
X^d)^d, X^d]^dY ,  X ] \overset {i} = 1$. Since $X \equiv B^{k_B}$, it
follows from Lemma \ref{L2B} that $[((X^dY^d )^d X^d)^d, X^d]^dY \overset {i}
= B^\ell$ for some $\ell$ and so $G^{k_G}$ is  conjugate in rank $i$ to
$B^\ell$. By Lemmas \ref{L2} and \ref{L4},
$$
| B^{k_B} | \le  \zeta^{-5}d^{-3}|F^{k_F}| \le
\zeta^{-6}d^{-4}|G^{k_G}| \le
10 \zeta^{-7}d^{-4} | G | < | G|
$$
(LPP: $\zeta \succ \eta = d^{-1}$), whence $|B| < |G|$. This, however,
contradicts  Lemma \ref{L2B}. Hence,  $[([((X^dY^d )^d X^d)^d,
X^d]^dY)^d,X^d] \overset {i} \neq 1$ and so $k_H \neq 0$.

By definitions, $ [ W_G  G^{d k_G } W_G^{-1},  B^{d k_B } ]  \overset {i}  =
W_H H^{k_H} W_H^{-1}$. By Lemma \ref{L2C}, $| k_H | \le 100 \zeta^{-1}$, $|
G^{d k_G} | \le \zeta^{-1} | H^{k_H} |$ and estimates (\ref{L5:1}) are
proved.

As in proofs of Lemmas \ref{L2}, \ref{L4}, we have from Lemmas \ref{L1} and
22.1 \cite{O89} that
$$
| W_H | <   ( \gamma + \tfrac 12 )\cdot 2( 2| W_G | +
 | G^{d k_G} | + | B^{d k_B } |  +  \tfrac 12  |  H^{k_H}  |) .
$$
Hence, by Lemmas \ref{L2}, \ref{L4} and estimates (\ref{L5:1}),
\begin{multline*}
| W_H | <   (1 + 2 \gamma ) (  (2 \zeta^{-1} d^{-1} +1 )| G^{d k_G} | +
\zeta^{-6} d^{-4}| G^{d k_G} | +   \tfrac 12  |  H^{k_H}  |)< \\ < (1 + 2
\gamma ) ((  2 \zeta^{-1} d^{-1} +1  + \zeta^{-6} d^{-4} )\zeta^{-1} + \tfrac
12 )  |  H^{k_H}  | <     2 \zeta^{-1}| H^{k_H} |
\end{multline*}
(LPP: $\gamma \succ \zeta \succ \eta = d^{-1}$), as required in (\ref{L5:2}).
Lemma \ref{L5} is proved.
\end{proof}

\begin{lemma}
\label{L6} Let $R_{A^f, j, z^*}$ be a defining word of rank $i+1$ defined by
$(\ref{R1})$ if $z^* = 1$ or by   $(\ref{R2})$ if $z^* = 2$. Then $0 < | f|
\le 100 \zeta^{-1}$, $|A| > d$, the words $T_{A^f, j, z^*}$, $U_{A^f, j,
z^*}$  do not belong to the cyclic subgroup $\la A \ra$, generated by $A$, of
$G(i)$ and
$$
\max(|T_{A^f, j, z^*} | , |U_{A^f, j, z^*} |) < d |A| .
$$
\end{lemma}

\begin{proof}    First we let $z^* = 1$. It follows from definitions that,
in the foregoing notation,  we can assume that
$$
A \equiv G ,   \quad
T_{A^f, j, 1}   \overset {i} = W_G^{-1} B^{k_B} W_G , \quad
U_{A^f, j, 1}   \overset {i} = W_G^{-1} Z  C^{k_C} Z^{-1} W_G ,
$$
and $f = f(A^f, j, 1)$ is $k_G$. Hence, in view of Lemmas \ref{L2},
\ref{L4},
\begin{gather*}
0 < | f| \le  100  \zeta^{-1} , \\
|A| \ge 10^{-1} \zeta^2 | F^{d k_F} | \ge 10^{-1} \zeta^7 d^3 | B^{d k_B} |
\ge
10^{-1} \zeta^7 d^3 > d
\end{gather*}
(LPP: $\zeta \succ \eta = d^{-1}$), and
\begin{multline*}
|T_{A^f, j, 1} | , |U_{A^f, j, 1} | \le 2( |W_G| +|Z|) + | B^{k_B}| +  |
C^{k_C}| < \\
< ( 2( \zeta^{-1} + 3 \zeta^{-6} d^{-3}) + 2 \zeta^{-6} d^{-4})  | G^{k_G} |
< 3 \zeta^{-1} | G^{k_G} | \le  30 \zeta^{-2} | G|< d |A|
\end{multline*}
(LPP: $\zeta \succ \eta = d^{-1}$).

Assume that one of $T_{A^f, j, 1}$, $U_{A^f, j, 1}$  belongs to
$\la G \ra \subseteq G(i)$. Then one of $B^{k_B}$, $C^{k_C}$
is  conjugate in rank $i$ to  a power of $G$. However, by Lemmas \ref{L2},
\ref{L4},
$$
|B^{k_B}|, |C^{k_C}| \le  \zeta^{-6} d^{-4} | G^{k_G} | \le 10 \zeta^{-7}
d^{-4} | G | < |G|
$$
(LPP: $\zeta \succ \eta = d^{-1}$), whence $|B|, |C| < |G|$ which is a
contradiction to Lemma \ref{L2B}.

Now we let $z^* = 2$. It follows from definitions that,
in the foregoing notation,  we can assume that
$$
A \equiv H,  \quad
T_{A^f, j, 2}   \overset {i} = W_H^{-1} W_G G^{k_G} W_G^{-1} W_H , \quad
U_{A^f, j, 2}   \overset {i} = W_H^{-1} Z  C^{k_C} Z^{-1} W_H ,
$$
and $f = f(A^f, j, 2)$ is $k_H$. Hence, in view of Lemmas \ref{L2}, \ref{L4},
\ref{L5},  we have
\begin{gather*}
0 < | f | \le 100  \zeta^{-1}  , \\
|A| \ge 10^{-2}  \zeta^2 |G^{d k_G} | \ge
  10^{-2}  \zeta^3 d |F^{d k_F} | \ge
10^{-2} \zeta^3 d^2 > d
\end{gather*}
(LPP: $\zeta \succ \eta = d^{-1}$) and
\begin{multline*}
|T_{A^f, j, 2} | , |U_{A^f, j, 2} | \le 2( |W_H| + |W_G| +|Z|) +
| G^{k_G}| +  | C^{k_C}| <    \\
< ( 2( 2\zeta^{-1}   +    \zeta^{-2} d^{-1}    +3 \zeta^{-7} d^{-4} )
+ \zeta^{-1} d^{-1} +  \zeta^{-7} d^{-5})  | H^{k_H} | <  \\
< 5 \zeta^{-1} | H^{k_H} | \le 500 \zeta^{-2} | H | <  d |A|
\end{multline*}
(LPP: $\zeta \succ \eta = d^{-1}$), as required.

Assume that one of $T_{A^f, j, 2}$, $U_{A^f, j, 2}$  belongs to
$\la H \ra \subseteq G(i)$. Then one of $G^{k_G}$, $C^{k_C}$
is  conjugate in rank $i$ to  a power of $H$. However, we saw above that
$|B|, |C| < |G|$ and it follows from Lemma \ref{L5} that
$$
|G| \le \zeta^{-1} d^{-1} \cdot 100 \zeta^{-1} |H| < |H|
$$
(LPP: $\zeta \succ \eta = d^{-1}$). Therefore, $|G|, |C| < |H|$ and, as
before, we have a contradiction to Lemma \ref{L2B}. Lemma \ref{L6} is proved.
\end{proof}

\begin{lemma} 
\label{L7} Suppose that $(X_1,Y_1)$ and $(X_2,Y_2)$ are two pairs of words
such that $(X_1^dY_1^d )^d X_1^d$ is conjugate in rank $i$ to $(X_2^dY_2^d
)^d X_2^d$ and $[X_1, Y_1] \overset {i} \neq 1$, $[X_2, Y_2] \overset {i}
\neq 1$. Then the pairs $(X_1,Y_1)$ and $(X_2,Y_2)$ are conjugate in rank
$i$.
\end{lemma}

\begin{proof}
Without loss of generality we can assume that
\begin{gather*}
X_1\equiv B_1^{k_{B_1}} , \quad  X_2\equiv B_2^{k_{B_2}} , \quad  \ Y_1
\equiv Z_1C_1^{k_{C_1}}Z_1^{-1} , \quad Y_2 \equiv
Z_2C_2^{k_{C_2}}Z_2^{-1} , \\
X_1^dY_1^d \equiv W_{D_1}D_1^{k_{D_1}}W_{D_1}^{-1} , \quad  X_2^dY_2^d \equiv
W_{D_2}D_2^{k_{D_2}}W_{D_2}^{-1} ,
\end{gather*}
where $B_1$, $B_2$, $C_1$, $C_2$, $D_1$, $D_2$ are words simple in rank
$\,i\,$ or periods of rank $\le i$, the words $Z_1$, $Z_2$, $W_{D_1}$,
$W_{D_2}$ have the minimal lengths among all words satisfying the
corresponding equalities, and if $A_1 \in \{ B_1, B_2, C_1, C_2, D_1, D_2\}$
is conjugate in rank $i$ to $A_2^{\pm 1}$, where $A_2 \in \{ B_1, B_2, C_1,
C_2, D_1, D_2\}$, then $A_1 \equiv A_2$.

Consider a reduced annular diagram $\Delta$ of rank $i$ for conjugacy of the
words $W_{D_1}D_1^{dk_{D_1}}W_{D_1}^{-1}B_1^{dk_{B_1}}$ and
$W_{D_2}D_2^{dk_{D_2}}W_{D_2}^{-1}B_2^{dk_{B_2}}$. Denote two cyclic sections
of the boundary $\p \Delta$ of $\D$ by $p_1q_1$ and $(p_2q_2)^{-1}$, where
\begin{gather*}
\ph(p_1) \equiv D^{d k_{D_1}}, \quad  \ph(q_1) \equiv
W_{D_1}^{-1}B_1^{dk_{B_1}}W_{D_1} , \\
\ph(p_2) \equiv D^{d k_{D_2}}, \quad \ph(q_2) \equiv
W_{D_2}^{-1}B_2^{dk_{B_2}}W_{D_2} .
\end{gather*}

It follows from  Lemma \ref{L2} that
\begin{gather*}
|W_{D_1}^{-1}B_1^{dk_{B_1}}W_{D_1}|<11\zeta^{-2}| D_1^{k_{D_1}}|\le
11\zeta^{-2} d^{-1}| D_1^{d k_{D_1}}| < \zeta^{3}| D_1^{d k_{D_1}}|
\end{gather*}
(LPP: $\zeta \succ \eta = d^{-1}$). Analogously,
$|W_{D_2}^{-1}B_2^{dk_{B_2}}W_{D_2}|< \zeta^{3}| D_2^{d k_{D_2}}|$. Hence,
\begin{equation}
|q_1| < \zeta^3|p_1|, \quad |q_2| < \zeta^3|p_2| . \label{L7:1}
\end{equation}

Assume that $|q_1|\ge \zeta |p_2|$. Then it follows from (\ref{L7:1}) that
$$
|p_2q_2|<(1+\zeta^3)|p_2|\le
\zeta^{-1}(1+\zeta^3)|q_1|<\zeta^2(1+\zeta^3)|p_1|<\zeta|p_1|.
$$
Hence, $|p_2q_2|<\zeta|p_1|$ and, in view of (\ref{L7:1}), $\Delta$ is an
$F$-map whose existence contradicts Lemma \ref{L1} and Lemmas 24.7, 25.10
\cite{O89}.

Therefore, we can assume that $|q_1|< \zeta |p_2|$ and, similarly, $|q_2|<
\zeta |p_1|$. Now we can see that  $\Delta$ is an $I$-map. It follows from
Lemmas \ref{L2I} and 25.10 \cite{O89} that $D_1\equiv D_2\equiv D$ and the
sections $p_1$, $p_2$ are $D$-compatible. By Lemma \ref{L2}, $|k_{D_1}|\le
100\zeta^{-1}$. Hence, using estimate (\ref{L7:1}), we get
$$
|q_1|<\zeta^3|p_1|=\zeta^3|D_1^{dk_{D_1}}|\le
100\zeta^{2}|D_1^{d}|<\zeta|D^{d}_1| =\zeta|D^{d}| .
$$
Analogously, $|q_2| <\zeta|D^{d}|$.

If $k_{D_1}\ne k_{D_2}$, then cutting $\Delta$ along a simple path, that
makes $p_1$ and $p_2$ $D$-compatible, we could turn $\Delta$ into an $F$-map
whose existence  contradicts Lemma \ref{L1}  and Lemmas 24.7, 25.10
\cite{O89}. Hence, it is shown that $k_{D_1}=k_{D_2}=k$. Cutting $\Delta$
along a simple path, that makes $p_1$ and $p_2$ $D$-compatible, we can see
that
$$
W_{D_2}D^{\ell_D} W_{D_1}^{-1}B_1^{dk_{B_1}}W_{D_1}D^{-\ell_D}W_{D_2}^{-1}
\overset i = B_2^{dk_{B_2}}
$$
for some integer $\ell_D$. By Lemma \ref{L2B}, $B_1$ and $(B_2)^{\pm 1}$ are
conjugate in rank $i$. Hence, $B_1\equiv B_2\equiv B$. It also follows from
Lemma \ref{L2B} that $k_{B_1}=k_{B_2}$ and $W_{D_2}D^{\ell_D}W_{D_1}^{-1}
\overset i = B^{\ell_B}$ for some integer $\ell_B$. This implies that
$W_{D_1} \overset i = B^{-\ell_B}W_{D_2}D^{\ell_D}$ and so
\begin{multline*}
X_1^dY_1^d \overset i = W_{D_1}D^k W_{D_1}^{-1} \overset i =
B^{-\ell_B}W_{D_2}D^{\ell_D} D^k D^{-\ell_D}W_{D_2}^{-1}B^{\ell_B} \overset i
= \\
\overset i =  B^{-\ell_B}X_2^dY_2^dB^{\ell_B} \overset i =
X_2^dB^{-\ell_B}Y_2^dB^{\ell_B}.
\end{multline*}
Since $X_1 \equiv  X_2 \equiv  B^{k_{B_1}}$, we obtain that  $Y_1^d \overset
i = B^{-\ell_B}Y_2^dB^{\ell_B}$. By Lemma \ref{L2A}, $Y_1 \overset i =
B^{-\ell_B} Y_2 B^{\ell_B}$ and we see that the pair $(X_2, Y_2)$ is
conjugate in rank $i$ to $(X_1, Y_1)$ by $B^{\ell_B}$. Lemma \ref{L7} is
proved.
\end{proof}

\begin{lemma}
\label{L9} Suppose that $(X_1,Y_1)$ and $(X_2,Y_2)$ are two pairs of words
such that $(X_1^dY_1^d )^d X_1^d$ is conjugate in rank $i$ to $((X_2^dY_2^d
)^d X_2^d)^{-1}$ and $[X_1, Y_1] \overset {i} \neq 1$, $[X_2, Y_2] \overset
{i} \neq 1$. Then the pairs $(X_1, Y_1)$ and $(X_2^{-1}, Y_2^{-1})$ are
conjugate in rank $i$.
\end{lemma}

\begin{proof} Observe that $((X_2^dY_2^d )^d
X_2^d)^{-1} \overset 0 = ((X_2^{-1})^d(Y_2^{-1})^d )^d (X_2^{-1})^d$.
Therefore, if $(X_1^dY_1^d )^d X_1^d$ is conjugate in rank $i$ to
$((X_2^dY_2^d )^d X_2^d)^{-1}$, then $(X_1^dY_1^d )^d X_1^d$ is conjugate in
rank $i$ to $((X_2^{-1})^d(Y_2^{-1})^d )^d (X_2^{-1})^d$. Hence, our claim
follows from Lemma~\ref{L7}.
\end{proof}

\begin{lemma} 
\label{L10} Suppose that $(X_1,Y_1)$ and $(X_2,Y_2)$ are two pairs of words
such that $v_1(X_1,Y_1)$ is conjugate in rank $i$ to $v_1(X_2,Y_2)$ and
$[X_1, Y_1] \overset {i} \neq 1$, $[X_2, Y_2] \overset {i} \neq 1$. Then the
pairs $(X_1, Y_1)$ and $(X_2, Y_2)$ are conjugate in rank $i$.
\end{lemma}

\begin{proof}
Conjugating the pairs $(X_1, Y_1)$, $(X_2, Y_2)$ in rank $i$ if necessary,
we can assume that
$$
X_1 \equiv B_1^{k_{B_1}} , \quad   Y_1 \equiv Z_1 C_1^{k_{C_1}} Z_1^{-1} ,
\quad  X_2 \equiv B_2^{k_{B_2}} ,  \quad  Y_2 \equiv Z_2 C_2^{k_{C_2}}
Z_2^{-1} ,
$$
where each of  $B_1$, $B_2$, $C_1$, $C_2$ is either simple in rank $i$ or a
period of some rank $\le i$ and, when $B_1^{k_{B_1}}$, $B_2^{k_{B_2}}$,
$C_1^{k_{C_1}}$, $C_2^{k_{C_2}}$ are fixed, the words $Z_1$, $Z_2$ are picked
to have minimal length. Furthermore, consider the following equalities
\begin{gather*}
X_1^dY_1^d \overset {i} = W_{D_1} D_1^{k_{D_1}} W_{D_1}^{-1} ,  \quad
X_2^dY_2^d \overset {i} = W_{D_2} D_2^{k_{D_2}} W_{D_2}^{-1} ,  \\
(X_1^dY_1^d )^d X_1^d \overset {i} = W_{E_1} E_1^{k_{E_1}} W_{E_1}^{-1} ,
\quad (X_2^dY_2^d )^d X_2^d \overset {i} = W_{E_2} E_2^{k_{E_2}} W_{E_2}^{-1}
,
\\ [((X_1^dY_1^d )^d X_1^d)^d,  X_1^d]  \overset {i} = W_{F_1}
F_1^{k_{F_1}} W_{F_1}^{-1} , \quad [((X_2^dY_2^d )^d X_2^d)^d,  X_2^d]
\overset {i} = W_{F_2} F_2^{k_{F_2}} W_{F_2}^{-1} ,
\end{gather*}
where each of $D_1$, $E_1$, $F_1$, $D_2$, $E_2$ $F_2$ is either a simple in
rank $i$ word or a period of some rank $\le i$ and the conjugating words
$W_{D_1}$, $W_{E_1}$, $W_{F_1}$, $W_{D_2}$, $W_{E_2}$, $W_{F_2}$ are picked
(when $D_1$, $E_1$, $F_1$, $D_2$, $E_2$ $F_2$ are fixed) to have minimal
lengths.

We can also assume that if $A_1 \in \{ B_1, C_1, D_1, E_1, F_1, B_2, C_2,
D_2, E_2, F_2\}$ is conjugate in rank $i$ to $A_2^{\pm 1}$, where $A_2 \in \{
B_1, C_1, D_1, E_1, F_1, B_2, C_2, D_2, E_2, F_2\}$, then $A_1 \equiv A_2$.

Consider a reduced annular diagram $\Delta$ of rank $i$ for conjugacy of the
words  $W_{F_1} F_1^{dk_{F_1}} W_{F_1}^{-1}Z_1 C_1^{k_{C_1}} Z_1^{-1}$ and
$W_{F_2} F_2^{dk_{F_2}} W_{F_2}^{-1}Z_2 C_2^{k_{C_2}} Z_2^{-1}$. Denote two
cyclic sections of the boundary $\p \D$ of $\D$ by $p_1q_1$ and
$(p_2q_2)^{-1}$,  where
\begin{gather*} \ph(p_1) \equiv F_1^{dk_{F_1}} , \quad
 \ph(q_1) \equiv W_{F_1}^{-1}Z_1 C_1^{k_{C_1}} Z_1^{-1}W_{F_1} , \\
\ph(p_2) \equiv F_2^{dk_{F_2}} , \quad \ph(q_2) \equiv  W_{F_2}^{-1}Z_2
C_2^{k_{C_2}} Z_2^{-1}W_{F_2} .
\end{gather*}
It follows from Lemma \ref{L2} that
\begin{multline*}
|W_{F_1}^{-1}Z_1 C_1^{k_{C_1}}
Z_1^{-1}W_{F_1}|<(4\zeta^{-1}+7\zeta^{-5}d^{-2})|F_1^{k_{F_1}}| \le \\ \le
(4\zeta^{-1}+7\zeta^{-5}d^{-2})d^{-1}|F_1^{d k_{F_1}}| < \zeta^3
|F_1^{dk_{F_1}}|
\end{multline*}
(LPP: $\zeta \succ \eta = d^{-1}$). Analogously, $|W_{F_2}^{-1}Z_2
C_2^{k_{C_2}} Z_2^{-1}W_{F_2}|< \zeta^3 |F_2^{dk_{F_2}}|$. Hence,
\begin{equation}
|q_1| < \zeta^3|p_1| , \quad |q_2| < \zeta^3|p_2|.  \label{L10:1}
\end{equation}

Assume that  $|q_1|\ge \zeta |p_2|$. Then it follows from (\ref{L10:1}) that
$$
|p_2q_2|<(1+\zeta^3)|p_2|\le
\zeta^{-1}(1+\zeta^3)|q_1|<\zeta^2(1+\zeta^3)|p_1|<\zeta|p_1|.
$$
Hence,  $|p_2q_2|<\zeta|p_1|$ and, in view of (\ref{L10:1}), $\Delta$ is an
$F$-map whose existence contradicts Lemma \ref{L1}  and Lemmas 24.7, 25.10
\cite{O89}.

Therefore, we can assume that $|q_1|< \zeta |p_2|$ and, similarly, $|q_2|<
\zeta |p_1|$. Now we can see that  $\Delta$ is an $I$-map. By Lemmas
\ref{L2I} and 25.10 \cite{O89}, $F_1\equiv F_2\equiv F$ and sections $p_1$
and $p_2$ are $F$-compatible. By Lemma \ref{L2}, $|k_{F_1}| \le
100\zeta^{-1}$. Hence, using estimate (\ref{L10:1}), we have
$$
|q_1|<\zeta^3|p_1|=\zeta^3|F_1^{dk_{F_1}}|\le
100\zeta^{2}|F_1^{d}|<\zeta|F^{d}_1| = \zeta|F^{d}|.
$$
Analogously, $|q_2|<\zeta|F^{d}|$.

If $k_{F_1}\ne k_{F_2}$, then, cutting $\Delta$ along a simple
path, that makes $p_1$ and $p_2$ $F$-compatible, we could turn
$\Delta$ into an $F$-map whose existence contradicts Lemma
\ref{L1} and Lemmas 24.7, 25.10 \cite{O89}. Hence,
$k_{F_1}=k_{F_2}$ and, cutting $\Delta$ along a simple path,
that makes $p_1$ and $p_2$ $F$-compatible,  we can see that
$Y_1$ and $Y_2$ are conjugate in rank $i$.

Since  $k_{F_1}=k_{F_2}$, it follows from definitions that the
word $[W_{E_1} E_1^{dk_{E_1}} W_{E_1}^{-1}, B_1^{dk_{B_1}}]$ is
conjugate in rank $i$ to $[W_{E_2} E_2^{dk_{E_2}} W_{E_2}^{-1},
B_2^{dk_{B_2}}]$. Let $\Delta$ be a reduced annular diagram of
rank $i$ for conjugacy of these two words. Denote two cyclic
sections of the boundary $\p \D$ of $\D$ by $p_1q_1p_2q_2$ and
$(p_3q_3p_4q_4)^{-1}$,  where
\begin{gather*}
\ph(p_1) \equiv \ph(p_2)^{-1} \equiv  E_1^{d k_{E_1}} , \quad
\ph(q_1) \equiv \ph(q_2)^{-1} \equiv  W_{E_1}^{-1} B_1^{dk_{B_1}} W_{E_1}  \\
\ph(p_3) \equiv \ph(p_4)^{-1} \equiv  E_2^{d k_{E_2}} , \quad
\ph(q_3) \equiv \ph(q_4)^{-1} \equiv  W_{E_2}^{-1}
B_2^{dk_{B_2}} W_{E_2}  .
\end{gather*}
It follows from Lemma \ref{L2} that
\begin{multline*}
|W_{E_1}^{-1} B_1^{dk_{B_1}} W_{E_1}|<(4\zeta^{-2}+\zeta^{-4}d^{-1} ) |
E_1^{k_{E_1}}| \le \\ \le (4\zeta^{-2}+\zeta^{-4}d^{-1} ) d^{-1} | E_1^{d
k_{E_1}}| <\zeta^3| E_1^{dk_{E_1}}|
\end{multline*} (LPP: $\zeta \succ \eta = d^{-1}$). Therefore,
\begin{equation}
|q_1|, |q_2| < \zeta^3|p_1| \quad \mbox{and} \quad |q_1|, |q_2| <
\zeta^3|p_2|. \label{L10:3}
\end{equation}
Similarly,
\begin{equation}
|q_3|, |q_4| < \zeta^3|p_3| \quad \mbox{and} \quad |q_3|, |q_4| <
\zeta^3|p_4|. \label{L10:4}
\end{equation}

Assume that $|q_1|\ge \zeta |p_3|$.  Then, by (\ref{L10:3})--(\ref{L10:4}),
we have
$$
|p_3q_3p_4q_4|<2(1+\zeta^3)|p_3|\le
2\zeta^{-1}(1+\zeta^3)|q_1|<2\zeta^2(1+\zeta^3)|p_1|<\zeta|p_1|.
$$
In view of estimates (\ref{L10:3})--(\ref{L10:4}), we can turn
$\Delta$ into an $E$-map $\Delta'$ by pasting together $q_1$
and $q_2$. It follows from Lemmas 24.6, 25.10 \cite{O89} that
the images of $p_1$ and $p_2$ in $\Delta'$ are
$E_1$-compatible. This implies that $k_{F_2}=0$, contrary to
Lemma \ref{L2}.

Therefore, we can assume that $|q_1|< \zeta |p_3|$ and,
similarly, $|q_3|< \zeta |p_1|$. Now we can see that $\Delta$
is an $I$-map. If  $p_1$ and $p_2$ are $E_1$-compatible in
$\Delta$, then $k_{F_2}=0$, contrary to Lemma \ref{L2}. Hence,
$p_1$ and $p_2$ may not be  $E_1$-compatible in $\Delta$.
Similarly, $p_3$ and $p_4$ are not $E_2$-compatible in
$\Delta$. Therefore, it follows from Lemmas \ref{L2I} and
25.10 \cite{O89} that $E_1\equiv E_2\equiv E$ and either $p_1$
is $E$-compatible with $p_3^{-1}$ and $p_2$ is $E$-compatible
with $p_4^{-1}$ or $p_1$ is $E$-compatible with $p_4^{-1}$ and
$p_2$ is $E$-compatible with $p_3^{-1}$. Let $t_1$ and $t_2$ be
simple disjoint paths in $\D$ that make corresponding pairs of
paths $p_1$, $p_2$, $p_3^{-1}$, $p_4^{-1}$ $E$-compatible. Let
us cut $\D$ along $t_1$ and $t_2$ and then paste the  two
resulting diagrams along the images of $q_3$ and $q_4^{-1}$
(recall that $\ph(q_3) \equiv \ph(q_4)^{-1}$). Let $\D_0$
denote the simply connected diagram of rank $i$ thus obtained
from $\D$. Observe that
\begin{equation}
\ph (\p \D_0) \overset i = E_1^{\ell_E} W_{E_1}^{-1}
B_1^{dk_{B_1}}W_{E_1} E_1^{-\ell_E} W_{E_1}^{-1}
B_1^{-dk_{B_1}}W_{E_1} \overset i = 1 , \label{L10:5}
\end{equation}
where either $\ell_E = d(k_{E_1}  - k_{E_2})$ in the case when
$p_1$ is $E_1$-compatible in $\Delta$ with $p_3^{-1}$ and $p_2$
is $E_1$-compatible with $p_4^{-1}$  or $\ell_E = d(k_{E_1} +
k_{E_2})$ in the case when $p_1$ is $E_1$-compatible in
$\Delta$ with $p_4^{-1}$ and $p_2$ is $E_1$-compatible with
$p_3^{-1}$. If $k_{E_1} \neq \pm k_{E_2}$ then it follows from
Lemma \ref{L2B}, equality (\ref{L10:5}) and definitions that
$B_1 \equiv E_1$. On the other hand, it follows from Lemma
\ref{L2} that
$$
| B_1| \le \zeta^{-4} d^{-2} |E_1^{k_{E_1}} | < 100 \zeta^{-5}
d^{-2} |E_1 | < |E_1 |
$$
(LPP: $\zeta \succ \eta = d^{-1}$). This contradiction shows
that $k_{E_1} = \pm k_{E_2}$ and so $E_1^{k_{E_1} } \equiv
E_2^{\pm k_{E_2}}$. By definitions, this means that the word
$(X_1^dY_1^d )^d X_1^d$ is conjugate in rank $i$ to
$((X_2^dY_2^d )^d X_2^d)^{\pm 1}$.

Assume that $(X_1^dY_1^d )^d X_1^d$ is conjugate in rank $i$ to
$((X_2^dY_2^d )^d X_2^d)^{-1}$.  Then, by Lemma \ref{L9}, $Y_1$
is conjugate in rank $i$ to $Y_2^{-1}$. On the other hand, as
we saw above, $Y_1$ is conjugate in rank $i$ to $Y_2$. Hence,
$Y_2$ is conjugate in rank $i$ to $Y_2^{-1}$. This, however, by
Lemmas \ref{L2B} and \ref{L2A}, implies that $Y_2 \overset i =
1$. This contradiction to $[X_1, Y_1] \overset {i} \neq 1$
proves that $(X_1^dY_1^d )^d X_1^d$ is conjugate in rank $i$ to
$(X_2^dY_2^d )^d X_2^d$. Now a reference to Lemma \ref{L7}
completes the proof of Lemma \ref{L10}.
\end{proof}

\section{Proof of Theorem}

Our Theorem is immediate from the following below Lemmas
\ref{L11}--\ref{L12}.

\begin{lemma} 
\label{L11} The group $G(\infty)$, defined by presentation
\eqref{G8}, is naturally isomorphic to the free group $F(\A) /
W_{1,2}(F(\A))$ of the variety $\MM$ in the alphabet $\A$.
\end{lemma}

\begin{proof}  It follows from the definition  of defining
words of the group $G(\infty)$ that each of them is in
$W_{1,2}(F(\A))$  and so there is a natural epimorphism
$$
G(\infty) \to F(\A) / W_{1,2}(F(\A)) .
$$

Suppose that $\wtl X$, $\wtl Y$ are some words in $\A^{\pm 1}$ and
\begin{equation}
w_{z^*}(\wtl X, \wtl Y )   \neq  1  \label{L11:1}
\end{equation}
in the group $G(\infty)$, where $z^* \in \{ 1,2\}$.

Let $A$ be a period of some rank such that $A^f$ for some $f$
is conjugate in $G(\infty)$ to $v_{z^*}(\wtl X, \wtl Y )$. (The
existence of such an $A$ follows from definitions and Lemma
\ref{L2B}; see also Lemma 18.1 \cite{O89}.)  Note that, in view
of (\ref{L11:1}), $[\wtl X, \wtl Y ]   \neq 1$ in $G(\infty)$.
Hence, by Lemmas \ref{L2}, \ref{L4}, \ref{L5}, we can replace
the pair $(\wtl X, \wtl Y )$ by a conjugate in the group
$G(\infty)$ pair $(X, Y)$ such that $X \equiv B^{k_B}$, $Y
\equiv Z C^{k_C} Z^{-1}$, and $v_{z^*}(X, Y ) = W_A A^f
W_A^{-1}$ in $G(\infty)$, where $B$, $C$ are some periods,
$|f|>0$ and
\begin{gather*}
|X^d | + |Y^d | = | B^{d k_B} | + |C^{d k_C} | +2 |Z| <
8\zeta^{-5}d^{-2} | F^{k_F} | \le  8 \zeta^{-6}d^{-3} |A^f| .
\end{gather*}
Hence, the length $ | v_{z^*}( X, Y ) | $, which is either
$|[((X^dY^d )^d X^d)^d,  X^d]^dY|$ if $z^* = 1$ or
\newline
$|[([((X^dY^d )^d X^d)^d,  X^d]^dY)^d,X^d]|$ if $z^* = 2$, can
be estimated as follows
\begin{multline*}
| v_{z^*}( X, Y ) |  \le 2d(d(2d(d( |X^d| + |Y^d|) + |X^d|) + 2|X^d|)+|Y|) +
2|X^d| < \\
< 5 d^4( |X^d| + |Y^d|) < 40  d \zeta^{-6} |A^f| < 10^4 d \zeta^{-7} |A|
\end{multline*}
for  $0 < |f| \le 100 \zeta^{-1} $ by Lemmas \ref{L4},
\ref{L5}. Consider a reduced annular diagram $\D$ of some rank
$i'$ for conjugacy of $v_{z^*}( X, Y )$  and $A^f$. By Lemmas
\ref{L1} and 22.1 \cite{O89}, $\D$ can be cut into a simply
connected diagram $\D_1$ along a simple path $t$ which connects
points on distinct components of $\p \D$ with $| t| < \gamma
|\p \D |$. Therefore,
$$
|\p \D_1 | < (1 + 2 \gamma) |\p \D | < (1 + 2 \gamma)(10^4 d \zeta^{-7}+
100\zeta^{-1} )|A|
< \tfrac 12 n |A|
$$
(LPP: $\gamma \succ \zeta \succ \eta = d^{-1} \succ \iota =
n^{-1}$). Then, by Lemmas \ref{L1}, 20.4 and 23.16 \cite{O89}
applied to $\D_1$, the diagram $\D_1$ contains no 2-cells of
rank $ \ge |A|$, whence $\D_1$, $\D$ are diagrams of rank $|A|
-1$. Since $A \in \X_{|A|}$, it follows from the construction
of defining words of rank $|A|$ that there will be a defining
word in $\Ss_{|A|}$ which guarantees that $w_{z^*}(X, Y)
\overset {|A|} = 1$.  A contradiction to the assumption
(\ref{L11:1}) proves that $G(\infty)$ is in $\MM$ and Lemma
\ref{L11} is proved.
\end{proof}

\begin{lemma} 
\label{L12} The group $G(\infty)$, defined by presentation
\eqref{G8}, is not hopfian.
\end{lemma}

\begin{proof} By Lemma \ref{L11},   $G(\infty)$ is
the free group of $\MM$ in $\A$ and, therefore, every map $a_1
\to U_1, \dots$, $a_m \to U_m$, where $U_1, \dots, U_m$ are
words in $\A^{\pm 1}$, extends to a homomorphism $G(\infty) \to
G(\infty)$.

Consider a   homomorphism $\psi_{\infty} : G(\infty) \to
G(\infty)$ defined by $\psi_{\infty}(a_j) = a_j$ if $j \neq 2$
and $\psi_{\infty}(a_2) = v_1(a_1, a_2)$. Observe that the
relation $w_2(a_1, a_2) = 1$ in $G(\infty)$ and the definition
(\ref{w2}) of the word $w_2(x, y)$ ensure that  $a_2 \in \la
a_1, v_1(a_1, a_2) \ra \subseteq G(\infty)$. Hence,
$\psi_{\infty}$ is an epimorphism. Assume that $\psi_{\infty}$
is an automorphism. Then it follows from the relation $w_1(a_1,
a_2) = 1$ in $G(\infty)$  that the word
\begin{multline*}
U \equiv a_1^{\e_1} a_2^{n}  a_1^{\e_2} a_2^{n+2}  \dots a_1^{\e_{h/2-2}}
a_2^{n+h-6}
a_1^{\e_{h/2-1}} a_2^{n+h-4}   a_1^{\e_{h/2}} a_2^{(n+h-2) + h/2} \\
a_1^{\e_1} a_2^{-(n+1)}  a_1^{\e_2} a_2^{-(n+3)}  \dots a_1^{\e_{h/2-1}}
a_2^{-(n+h-3)} a_1^{\e_{h/2}} a_2^{-(n+h-1)}
\end{multline*}
is equal to 1 in $G(\infty)$.  Note that $U$ is cyclically
reduced and $|U| <(n+h)h$. Consequently,
\begin{equation}
|U| < (n+h) h < (1-\alpha)(h-1) n d \label{T:3}
\end{equation}
(LPP: $ \alpha \succ \zeta = h^{-1} \succ \eta = d^{-1} \succ \iota =
n^{-1}$).

On the other hand, since $U=1$ in $G(\infty)$ and $U \neq 1$ in
the free group $F(\A)$, it follows from Lemmas \ref{L1} and
23.16  \cite{O89} that
\begin{equation}
|U| >(1-\alpha) |\p \Pi | , \label{T:1}
\end{equation}
where $\Pi$ is a 2-cell of positive rank. Furthermore, it
follows from Lemmas \ref{L6} and 20.4 \cite{O89} (together with
condition $B$, see Sect. 20.4 \cite{O89}) that
\begin{equation}
|\p \Pi | > (h-1)n r(\Pi) > (h-1)n d \label{T:2}
\end{equation}
for an arbitrary 2-cell $\Pi$ of rank $r(\Pi) \ge 1$. However,
inequalities (\ref{T:1})--(\ref{T:2}) contradict (\ref{T:3}). This
contradiction shows that $U \neq 1$ in $G(\infty)$.  Thus, $U$ is
in the kernel of $\psi_{\infty}$ and $G(\infty)$ is not hopfian,
as required. \end{proof}

\smallskip

In conclusion, we remark that, for every $i \ge 0$, the group
$G(i)$, given by presentation (\ref{Gi}), is finitely presented
(this follows from definitions and Lemmas \ref{L2}--\ref{L6})
and satisfies a linear isoperimetric inequality (this can be
proved similar to Lemma 21.1 \cite{I94}). Therefore, $G(i)$ is
a torsion-free Gromov hyperbolic group (see \cite{G87}) and so,
by  Sela's  results \cite{Sl97}, $G(i)$ is a hopfian group.
Therefore, the non-hopfian group $G(\infty)$, given by
presentation (\ref{G8}), is a limit of hopfian groups $G(i)$,
$i=0,1,\dots$.

\end{document}